\documentclass[10pt]{article}
\usepackage{amsmath,amsthm,amsfonts,amssymb}
\usepackage{graphicx}
\usepackage{hyperref}
\usepackage[usenames]{color}

\newtheorem{theorem}{Theorem}[section]
\newtheorem{lemma}[theorem]{Lemma}
\newtheorem{claim}[theorem]{Claim}
\newtheorem{fact}[theorem]{Fact}
\theoremstyle{definition}

\newtheorem{corollary}[theorem]{Corollary}

\newtheorem{remark}[theorem]{Remark}

\newcommand{\comment}[1]{}

\def\wr{{\mathrm{wr}}}
\def\supp{{\mathrm{supp}}}

\title{Polynomial Time Conjugacy in Wreath Products and Free Solvable Groups}

\author{S. Vassileva}

\date{\today}

\begin{document}
\maketitle

\begin{abstract}

We prove that the complexity of the Conjugacy Problems for wreath products and for free solvable groups is decidable in polynomial time. 
For the wreath product $A \wr B$, we must assume the decidability in polynomial time of the Conjugacy Problems for $A$ and $B$ and of the power problem in $B$. We obtain the result by making the algorithm for the Conjugacy Problem in Matthews \cite{Matthews:1966} run in polynomial time. 
Using this result and properties of the Magnus embedding, we show that the Conjugacy and Conjugacy Search Problems in free solvable groups are computable in polynomial time.

\end{abstract}

\tableofcontents

\section{Introduction}


In this paper we attack the computational complexity of the
Conjugacy and Conjugacy Search Problems in free solvable groups. We
show that they are both solvable in polynomial time and that the degree of
the polynomial is uniform for the class of free solvable groups.
Further, we show that the Conjugacy Problem and the Conjugacy Search
Problem in wreath products are solvable in polynomial time modulo
some natural conditions.


Algorithmic problems in group theory were considered as early as
1910, when Dehn introduced the now famous Word and Conjugacy
Problems. Briefly, for a finitely generated group $G$, given two words as a product of generators, the Word Problem asks whether they are equal as elements of $G$ and the Conjugacy Problem asks whether they are conjugate to each other in $G$. Both of these decision problems quickly became an active area of research. Novikov (\cite{Novikov:1952}, \cite{Novikov:1958}) gave the first example of a finitely presented group with undecidable Word (and hence Conjugacy) Problem. A beautiful result of Miller exhibits a group which has decidable Word Problem and undecidable Conjugacy Problem \cite{Miller1}. At present, there are many interesting classes of groups where these problems are decidable. Here we mention only a few positive results about non-solvable groups and discuss  solvable groups in more detail below. The Word and Conjugacy  Problems are decidable in braid groups (Artin, \cite{Artin}), hyperbolic groups (Gromov, \cite{Gromov_hyperbolic}), wreath products of groups under some natural additional conditions (Matthews \cite{Matthews:1966}), the
Grigorchuk group (Grigorchuk \cite{Grig98}, Leonov \cite{Leonov}), bi-automatic groups (Gersten and Short, \cite{GerShort}), toral relatively hyperbolic groups,  free solvable groups
(Remeslennikov, Sokolov \cite{Remeslennikov-Sokolov:1970}).

Nowadays, while decidability is still an open area of research, the emphasis
has shifted to complexity of decidable problems. It is worth mentioning the work of
Lysenok, Miasnikov, and Ushakov who showed in \cite{LMU} that the
Conjugacy Problem in the Grigorchuk group is decidable in polynomial
time, the work the work of Lipton and Zalenstein on the polynomial time decidability
of the Word Problem in linear groups \cite{LipZal}, the work of Marshall,
Bridson and Haefliger, Epstein and Holt which, through successively
improving time bounds, culminates in showing that the Conjugacy
Problem in word-hyperbolic groups is decidable in linear time
\cite{Epstein_linhyp} and the work of Cannon, Goodman and Shapiro, and Holt and Rees \cite{Holt} in giving a linear time algorithm for deciding the Word Problem in nilpotent groups.


Solvable groups offer a whole new world on their own. An example of
Kharlampovich of a solvable group with undecidable Word, and hence
Conjugacy, Problem shows that one cannot derive any positive results
about the entire class of solvable groups. However, there are many
interesting subclasses in which the Conjugacy Problem is decidable,
for instance finitely generated metabelian groups (Noskov \cite{Noskov:1982}), nilpotent groups (Blackburn \cite{Blackburn_nilpotent}), polycyclic groups (Remeslennikov \cite{Remeslennikov_CP_polycyclic}) and free solvable groups
(Remeslennikov - Sokolov \cite{Remeslennikov-Sokolov:1970}. In all
of the above cases, however, the results are about decidability
without mention of the time complexity. The complexity of
algorithmic problems in solvable groups has recently become an active area of
research with a paper by Miasnikov, Roman'kov, Ushakov and Vershik
\cite{MRUV} which presents a cubic time algorithm to decide
the Word Problem in free solvable groups.


Most complexity results concern a fixed group. To the knowledge
of the author, there is no other studied class of infinite solvable groups for which
the Word and Conjugacy Problems can  be decided uniformly in
polynomial time. Even in the cases where one can solve the given
problem using a general description of the group, the algorithm
involves heavy pre-computations specific to this group which cannot
be generalized to produce a uniformly polynomial-time algorithm.


In this paper we use this result in \cite{MRUV} to show that the
Conjugacy Problem in free solvable groups is decidable in quintic
time. The proof follows the ideas of Remeslennikov and Sokolov
(\cite{Remeslennikov-Sokolov:1970}). First, we embed the free solvable
group of degree $(d+1)$ and rank $r$ in a wreath product of an
abelian group and a free solvable group of degree $d$. The image
of a word of length $n$ can be found in time $O(rdn^3)$. Since the images of two words under the Magnus embedding are conjugate if and only if these words are conjugate, we can apply our general result, namely that the Conjugacy Problem in this wreath product is decidable in polynomial time, provided the Conjugacy Problems in each factor (and the Power Problem in the second factor) are decidable in polynomial time. The second factor is a free solvable group of lesser degree, so we proceed by induction. Similarly, we solve the Conjugacy Search Problem. 

\comment{

The first section gives some general definitions and discusses known
algorithmic results on the Magnus embedding, which will be used
later on. In Section~\ref{sec:algorithm wreath} the algorithm given
by Matthews (\cite{Matthews:1966}) is modified and its complexity is
analyzed. In Section~\ref{sec:algorithm free solvable} it is shown
that the Conjugacy Problem in a free solvable group of degree $d$
and rank $r$ can be decided in time $O(rdn^3)$, where $n$ is the
length of the input words.
} 

\section{Preliminaries}

\subsection{Wreath products and the Magnus embedding}

We start by defining the objects essential to this paper -- wreath
products and the Magnus embedding.


Let $G$ be a group generated by a fixed finite set of generators $Y$. We
represent elements in $G$ by words $w$ over $Y^{\pm}$ and denote by $|w|$ the length of the word $w$.


Let $A$ and $B$ be groups. The \emph{restricted wreath
product} $A\wr B$ is the group formed by the set
$$A\wr B = \{bf \mid b\in B,\; f\in A^{(B)}\}, $$
with multiplication defined by $bfcg = bcf^cg$, where $f^c(x) = f(x
c^{-1})$ for $x \in B$, where $A^{(B)}$ denotes the set of functions
from $B$ to $A$ with finite support (i.e., functions from $B$ to $A$ which take non-zero values only for finitely many elements of $B$). Note that $A^{(B)}$ is a group under pointwise multiplication of functions with identity $1: B \rightarrow 1$, so we can view $A \wr B$ as the semi-direct product $B \ltimes A^{(B)}$.

Let $X = \{x_1, \hdots, x_n\}$ and $Y = \{y_1, \hdots, y_m\}$ be the
generating sets for $A$ and $B$, respectively. $A\wr B$ is generated
by $X, Y$ in the following sense: every function, $f \in A^{(B)}$
can be written as a product $f = \prod_i a_i^{b_i}$. Indeed consider
the functions of the form
$$f_{a_i,b_i}(x)=\left\{
\begin{array}{rl}
a_i & \text{if } x = b_i \\
1 & \text{otherwise }\\
\end{array} \right.$$
For simplicity, we denote $f_{a_i, 1}$ by $f_{a_i}$. Then for any $f
\in A^{(B)}$, one can write $f = \prod_i f_{a_i,b_i} = \prod_i
f_{a_i}^{b_i}$. There is clearly an identification between $f_{a_i}$
and $a_i$.

\begin{remark}
\label{remark: ordering supp(f)} One can rewrite a word $w = b_1
a_1\hdots b_k a_k$ in generators $X$ and $Y$ as $w = bf$ in
polynomial time. Observe that
$$w = b_1\hdots b_k a_1^{b_2\hdots b_k}\hdots a_2^{b_3\hdots b_k} a_{k-1}^{b_k} a_k.$$
Here $b = b_1\hdots b_k \in B$ and $a_1^{b_2\hdots b_k}
a_2^{b_3\hdots b_k}\hdots a_{k-1}^{b_k} a_k$ corresponds to a function in $A^{(B)}$ as follows. Denote $B_i = b_i\hdots b_k$. For each $1<i<j\leq k$, check whether
$B_i=B_j$ in $B$. This amounts to solving ${k-1}\choose 2$ Word Problems in
$B$. For each $B_{i_1} = B_{i_2} = \hdots = B_{i_j}$, write
$f(B_{i_1}) = a_1\hdots a_j$. This determines $f$ completely and we
can change presentations in time $O(|w|^2 T_{WB}(|w|))$, where
$T_{WB}$ is the time function for the Word Problem in $B$. Note that
if a word is given as a product of generators, converting it to
standard (or pair) form gives an ordering for $\supp(f) = \{B_i\}_i$
determined by the indices $i$. More precisely, $B_i < B_j$ whenever
$i < j$.
\end{remark}


Fix a free group $F$ of rank $r$ with basis $X$. The derived subgroup $F^{(d)}$ is defined by induction as follows: $F' = [F,F]$ and $F^{(d+1)} = [F^{(d)}, F^{(d)}]$. Define the free
solvable group, $S_{d,r} = F/F^{(d+1)}$.


Let $N$ be a normal subgroup of $F$. Denote by $\mu : F \rightarrow
F/N$ the canonical epimorphism. Let $U$ be a free
$\mathbb{Z}(F/N)$-module with basis $\{u_1, \hdots, u_r\}$, so $U
\simeq \mathbb{Z}(F/N) \oplus \hdots \oplus \mathbb{Z}(F/N) $. Then
the set of matrices

\begin{equation*}
M(F/N) = \left(
\begin{array}{cc}
F/N & U\\
0 & 1
\end{array} \right) = \bigg\{
\left(
\begin{array}{cc}
g & u\\
0 & 1
\end{array} \right)
\mid g\in F/N, u\in U \bigg\}
\end{equation*}
forms a group with respect to matrix multiplication. One can see
that (see for example, \cite{Remeslennikov-Sokolov:1970}) $M(F/N)
\simeq F/F' \wr F/N$.

The map $\varphi: F(X)\rightarrow M(F/N)$ defined by

$$x_i \mapsto \left(
\begin{array}{cc}
\mu (x_i) & u_i\\
0 & 1
\end{array}\right),  i=1, \hdots, r$$
extends to an injective homomorphism $\varphi : F/N^{\prime}
\rightarrow M(F/N)$, called the \emph{Magnus embedding}.

In the sequel, for $x \in F$ put

$$\varphi(x) =
\left(
\begin{array}{cc}
\mu (x) & u_x\\
0 & 1
\end{array} \right).$$

\subsection{Algorithmic Results for the Magnus Embedding}

Here we present and prove a few preliminary results on the Magnus
embedding that we will need in Section \ref{sec:algorithm free
solvable}.

\begin{theorem}[\cite{Remeslennikov-Sokolov:1970}]
\label{prop:f,g conj in S_(d+1) iff in M(S_d)} Let $\bar{f},\bar{g}
\in F/N'$, where $N$ is normal in $F$ and $N'$ is torsion-free. Then
$\bar{f}$ and $\bar{g}$ are conjugate in $F/N'$ if and only if their images in
$M(F/N)$ are conjugate.
\end{theorem}

In particular, the theorem above holds for the free solvable group
$F/F^{(d+1)}$, which is $F/N'$ for $N = F^{(d)}$.

\begin{theorem}[\cite{MRUV}]
\label{prop:magnus embedding is O(dr w^3)}
\label{prop:WP in S_d is O(dr w^3)}

The following hold:

\begin{enumerate}
\item[1)] For a given $w \in S_{d,r}$, one can compute $\varphi(w)$ in time
$O(dr|w|^3)$;
\item[2)] The Word Problem in $S_{d,r}$ is solvable in time $O(dr|w|^3)$, where $w$ is the input word.
\end{enumerate}
\end{theorem}

\begin{corollary}
\label{cor: reduction to CP in wreath is poly}
The Conjugacy Problem in $S_{d,r}$ reduces to the Conjugacy Problem in $F/F' \wr S_{d-1,r}$ it time $O(rdL^3)$, where $L$ is the length of the input words. 
\end{corollary}

The Power Problem for a group $G$ for given elements $x,y \in G$ consists of determining whether there exists an integer $n\in\mathbb{Z}$ such that $x = y^n$ and if so, to find it. 


\begin{theorem}
\label{prop:power problem is O(dr w^3)}

The power problem in $F/F^{(d)}$ is decidable in time $O(rdL^6)$,
where $r$ is the rank of $F$ and $L=|x|+|y|$ is the length of the input.
\end{theorem}

\begin{proof}  

Let $x$ and $y$ be elements in $F/F^{(d)}$ given as products of generators. Consider first the two trivial cases. If $y=1$, which can be checked in time $O(rd|y|^3)$, the problem reduces to a Word Problem, which is decidable in $O(rd|x|^3)$. If $x=1$, then $n=0$ is always a solution. Hence, after some preliminary computation which can be done in $O(rdL^3)$, we can assume without loss of generality that both $x$ and $y$ are non-trivial elements in $F/F^{(d)}$. Observe the following. 

\begin{fact}
\begin{enumerate}
\item \label{fact1}If there exists $n\in \mathbb{Z}$ such that $x=y^n$ in $F/F^{(d)}$, then $x=y^n$ in $F/F'$. 
\item \label{fact2}If there exists $n\in \mathbb{Z}$ such that $x=y^n$ in $F/F^{(d)}$, then $n$ is unique with this property. 
\end{enumerate}
\end{fact}

The first claim follows easily since $F/F'$ is a quotient of $F/F^{(d)}$ and the second one follows from the fact that free solvable groups are torsion-free. We proceed to solve the general case of the Power Problem in a free solvable group $F/F^{(d)}$. 

\begin{enumerate}
\item[Step 1: ] Solve the Power Problem in $F/F'$. It is a free
abelian group, so the elements $x$ and $y$ can be uniquely presented
in the form $x=x_1^{a_1} \hdots x_r^{a_r}$ and $y=x_1^{b_1} \hdots
x_r^{b_r}$, where $X = \{x_1, \hdots, x_r\}$ is the basis for $F$.
Obviously, this decomposition can be found in log-linear time, which is certainly in $O(rL^6)$. Then for each $1 \leq i \leq r$ set $n_i = a_i/b_i$. If all $n_i$ are equal
and integer, then $x=y^{n_1}$, as required. Otherwise, $x \not\in
\langle y\rangle$ and we are done. Clearly, this can be done in time
$O(r(|x|+|y|))$.

Note that the exponent $n$ satisfies $n \leq |x|+|y| = L$. 

\item[Step 2: ] Using $n$ from Step~1, check whether the equation 
\begin{equation}
\label{eqn: x=y^n}
x=y^n
\end{equation} 
holds in $F/F^{(d)}$. By Theorem~\ref{prop:WP in S_d is O(dr w^3)}, this can be done in time $O\big(rd \big( |x| + n|y|\big)^3 \big) \subseteq O( rd L^6)$. If (\ref{eqn: x=y^n}) does not hold, then $x \neq y^m$ for all integers $m$. Indeed, if there were some $m\in\mathbb{Z}$ for which $x=y^m$ in $F/F^{(d)}$, then by Fact\ref{fact1} the same equation would hold in $F/F'$. But by the uniqueness of $n$ (Fact\ref{fact2}), this is impossible. 
\end{enumerate}

\end{proof}

\section{Complexity of the Conjugacy Problem in Wreath Products}
\label{sec:algorithm wreath}

We establish a bound on the complexity of the Conjugacy Problem in
wreath products $A\wr B$ by giving a bound for a variant of the
algorithm developed by Matthews \cite{Matthews:1966}.

Let $x=bf, y=cg \in A\wr B$, where $b,c \in B$ and $f,g \in A$.
Denote $\supp(f) = \{b_1, \hdots, b_n\}$ and $\supp(g) = \{\beta_1,
\hdots, \beta_m\}$ where the $b_i$ and $\beta_j$ are ordered as in Remark ~\ref{remark: ordering supp(f)}. Recall that all elements are given as words in generators. Let $\bar{b}$ and $\bar{\beta}$ be the longest elements in $\supp(f)$ and in $\supp(g)$, and $\bar{a}$
and $\bar{\alpha}$ be the longest element in the image of $f$ and of $g$, respectively.

For each left $\langle b \rangle$-coset in $B$ that intersects $\supp(f)\cup \supp(g)$, choose a coset representative from $\supp(f) \cup \supp(g)$ and let $T_b=\{t_i\}_{i \in I_1\cup I_2}$, where $I_1$ indexes the coset representatives we just chose and $I_2$ indexes the remaining ones. Deciding whether $b_i, b_j \in \supp(f) \cup \supp(g)$ are in the
same coset is a Power Problem, since $b_i, b_j$ are in the same coset if and only if $b_i b_j^{-1} = b^k$ for some $k$. To find $T_b$ one needs to solve the Power Problem ${(n+m) \choose 2}$ times (for all pairs $(b_i, b_j)$). Hence it takes time ${(n+m) \choose 2}T_{PB}( 2|\bar{b}| + 2|\bar{\beta}| + |b|)$, where $T_{PB}$ is the time function for the power problem in $B$. For each $\gamma\in B$ and $i \in I_1 \cup I_2$, associate with $T_b$ the following map $\pi_{t_i}^{(\gamma)}: A^{(B)} \rightarrow A$:

\begin{equation*}
\pi_{t_i}^{(\gamma)}(f) = \left\{
\begin{array}{rl}
\prod\limits_{j=0}^{N-1}f(t_ib^j\gamma^{-1}) & \text{if } b \text{ is of finite order } N, \\
\\
\prod\limits_{j=-\infty}^{\infty}f(t_ib^j\gamma^{-1}) & \text{if } b \text{ is of infinite order. }\\
\end{array} \right.
\end{equation*}


Note that in the above all the products are finite, since $f$ has finite support. Denote $\pi_{t_i}^{(1)}(f)$ by $\pi_{t_i}(f)$. Matthews gives a condition to check conjugacy, which will be used here.

\begin{theorem}[\cite{Matthews:1966}]
\label{thm: conjugacy criterion for wreath}
Let $A$, $B$ be finitely
generated groups. Two elements $x=bf, y=cg \in A\wr B $ are
conjugate if and only if there exists $d\in B$ such that for all
$t_i \in T_b$ the following hold:
\begin{enumerate}
\item[(1)] $db = cd$,
\item[(2)] when the order of $b$ is finite, $\pi_{t_i}^{(d)}(g)$ is conjugate to $\pi_{t_i}(f)$ in A,
\item[(3)] when the order of $b$ is infinite, $\pi_{t_i}^{(d)}(g)= \pi_{t_i}(f)$ in A.
\end{enumerate}
\end{theorem}

In order to use this criterion computationally, we need to circumvent the use of the conjugator $d$. 

\begin{lemma}
\label{lemma: conjugacy bar and tilde}
Let $\{\bar{s_i}\}_{i\in I}$ and $\{\tilde{s_i}\}_{i\in I}$ be two sets of left $\langle c \rangle$-coset representatives such that $\bar{s_i} \langle c \rangle = \tilde{s_i} \langle c \rangle$. Then $\pi_{\bar{s_i}}(g)$ and $\pi_{\tilde{s_i}}(g)$ are conjugate for any $i\in I$. 
\end{lemma}

\begin{proof}
Since $\bar{s_i} \langle c \rangle = \tilde{s_i} \langle c \rangle$, there is some integer $k_i$ for which $\bar{s_i} = \tilde{s_i}c^{k_i}$ and hence, 
$$\pi_{\bar{s_i}}(g) = \prod_j g(\bar{s_i}c^j) = \prod_j g(\tilde{s_i}c^{k_i}c^j) = \prod_j g(\tilde{s_i}c^{k_i+j}).$$
This last product is a cyclic permutation of the factors in $\prod_j g(\tilde{s_i}c^j) = \pi_{\tilde{s_i}}(g)$ and so is conjugate to $\pi_{\tilde{s_i}}(g)$.
\end{proof}

Using the Theorem~\ref{thm: conjugacy criterion for wreath} and Lemma~\ref{lemma: conjugacy bar and tilde} we show that the time complexity of the Conjugacy Problem in wreath products is polynomial.

\begin{theorem}\label{thm: CP in wreath products is poly}
Let $A$ and $B$ be finitely generated groups such that the following
hold:
 \begin{enumerate}
 \item[1)] there are decision algorithms for the Conjugacy Problem in $A$ and in $B$ with polynomial time functions,
 $T_{CA}$, $T_{CB}$, respectively;
 \item[2)] there is an algorithm with polynomial time function $T_{PB}$ for the Power Problem in $B$.
 \end{enumerate}
 Then the Conjugacy Problem in $A\wr B$ is decidable with complexity

\begin{equation}
\label{eqn: complexity wreath prod}
O\big(L^2T_{CA}(L^2) + LT_{CB}(L) + L^2T_{PB}(L) \big), 
\end{equation}
where $L=|x| +|y|$ is the length of the input pair $x,y \in A\wr B$.
\end{theorem}

\begin{remark}
\label{remark: CP => WP} Note that every Word Problem ''s $x=1$?" is precisely the Conjugacy Problem ''Is $x$ conjugate to $1$"? To simplify the presentation, the complexities of
all Word Problems considered in this section will be bounded by the
complexities of the corresponding Conjugacy Problems.
\end{remark}

\begin{proof}
Let $x=bf, y=cg \in A\wr B$. The notation from the beginning of this
section will be used throughout. In order to simplify the subsequent treatment of complexity in this section, we will implicitly use the bounds 
$$|x|, |y|, n, m, |c|, |b|, |\bar{b}|, |t_i|, |\bar{a}| \leq L. $$

\begin{claim}
\label{subsec: compute pi_i} There is a polynomial time algorithm
which computes $\pi_{t_i}^{(\gamma)}(f)$. More precisely, 
\begin{itemize}
\item $\pi_{t_i}^{(\gamma)}(f)$ can be computed in time $L T_{PB}(L)$. 
\item $|\pi_{t_i}^{(\gamma)}(f)| \leq L^2$. 
\end{itemize}
\end{claim}

\begin{proof}

The algorithm is as follows:

\begin{description}
\item[Step 1: ] For each $b_k \in \supp(f)$ check whether there is some $j$ such that $t_i b^j \gamma ^{-1} = b_k$, i.e., $t_i^{-1}b_k \gamma = b^j$. This is an instance of the
Power Problem in $B$ and so can be done in time $T_{PB}(2|\bar{b}| + |b| + |\gamma|)$. If such $j$ exists, look up the corresponding value $a_j = f(b_k)$. Otherwise, $a_j$ does not occur
in the product.
\item[Step 2: ] There are $n$ elements in $\supp (f)$ to perform computations on,
so computing $\pi_{t_i}^{(\gamma)}(f)$ takes time $n T_{PB}(2|\bar{b}| + |b| + |\gamma|)$.
\item[Step 3: ] Set $\pi_{t_i}^{(\gamma)} = \prod_j a_j$. Note that the order in which the factors are multiplied is a priori determined by the solution $j$ to the Power Problem. However, if the order of $b$ is finite, by the definition of $\pi$ we take $j \mod N$, and if the order of $b$ is infinite, then the solution to the Power Problem is unique because in this case $b$ has no torsion. Thus, a fortiori, $\pi_{t_i}^{(\gamma)}$ is indeed equal to $\prod_j a_j$, where the $a_j$ are computed as above. 
\end{description}

Note that $|\pi_{t_i}^{(\gamma)}(f)| \leq n|\bar{a}|$, since each factor
in the product $\pi_{t_i}^{(\gamma)}(f)$ is in the image of $f$.

\end{proof}

\;\;\;

We modify the algorithm from \cite{Matthews:1966} so that it runs in polynomial time as follows:

\begin{description}
\item[Step 1.] Determine whether $b$ and $c$ are conjugate in
$B$. This takes time $T_{CB}(|x|+|y|) \in O(T_{CB}(L))$. If not, $x$
and $y$ are not conjugate. If $b$ and $c$ are conjugate in $B$, let
$d\in B$ be such that $db = cd$ (it is not required to find this
$d$).

\item[Step 2.] Consider the following three cases. 
\end{description}

\begin{description}
\item[Case 1:] $g = 1$. Then $\pi_{t_i}^{(d)}(g) = 1$, so $x$ and $y$ are
conjugate if and only if $\pi_{t_i}(f) = 1$. To check this compute
$\pi_{t_i}(f)$ as in Claim~\ref{subsec: compute pi_i} and solve the Word
Problem in $A$. This takes time
\begin{equation}\label{eqn: complexity case 1}
O \big( L T_{PB}(L) + T_{CA}(L^2) \big).
\end{equation}

\item[Case 2:] $g \neq 1$, and $\pi_{t_i}(f) = 1$ for all $i\in I_1$. In order
to check the latter, simply compute $\pi_{t_i}(f)$ for all $i\in I_1$. This
will take time $O(L^2 T_{PB}(L))$. Then, by Theorem~\ref{thm: conjugacy criterion for wreath}, $x$ is conjugate to $y$ if and only if $\pi_{t_i}^{(d)}(g) = 1$ for all $i\in I_1$ (since the $\pi_{t_i}^{(d)}(g) = 1$ for $i\in I_2$). Note that we need not know what $d$ actually is -- its existence is enough. Indeed, since $db=cd$, $g(t_i b^j d^{-1}) = g(t_i d^{-1}c^j)$ and hence 
$$\pi_{t_i}^{(d)}(g) = \prod_j g(t_i b^j d^{-1}) = \prod_j g(t_i d^{-1}c^j) = \pi_{t_id^{-1}}(g),$$
where $\{t_id^{-1}\}_{i\in I_1 \cup I_2}$ is a set of left $\langle c \rangle$-coset representatives. Moreover, by Lemma~\ref{lemma: conjugacy bar and tilde}, $\pi_{t_id^{-1}}(g)$ is conjugate to $\pi_{s_i}(g)$ for any other set of left $\langle c \rangle$-coset representatives $\{s_i\}_{i\in I_1 \cup I_2}$ for which $t_id^{-1}\langle c \rangle = s_i\langle c \rangle$. It follows that $\pi^{(d)}_{t_i}(g) = 1$ for all $i\in I_1 \cup I_2$ if and only if $\pi_{s_i}(g) = 1$ for all $i \in I_1\cup I_2$. 

Since $\pi_{s_i}(g)=1$ for all $i\in I_2$, to check whether $x$ and $y$ are conjugate, it is enough to check whether for some set of left $\langle c \rangle$-coset representatives $T_c = \{s_i\}_{i\in I_1}$, $\pi_{s_i}(g) = 1$ for all $i \in I_1$. Choosing $T_c$ can be done in time $O(L^2 T_{PB}(L))$ and by Claim~\ref{subsec: compute pi_i}, checking whether $\pi_{s_i}(g) = 1$ for all $i \in I_1$ can be done in time $L^2T_{CA}(L^2)$. Thus checking whether $x$ and $y$ are conjugate takes time

\begin{equation}\label{eqn: complexity case 2}
O\big( L^2 T_{PB}(L) + L^2 T_{CA}(L^2) \big).
\end{equation}

\item[Case 3:] $g \neq 1$ and some $\pi_{t_i}(f) \neq 1$. There are two
subcases:

\item[1)] \emph{The order of $b$ is finite.} By Theorem~\ref{thm: conjugacy criterion for wreath}, $x$ and $y$ are conjugate if and only if $\pi_{t_i}(f)$ and $\pi_{t_i}^{(d)}(g)$ are conjugate. As in Case~2, $\pi_{t_i}^{(d)}(g) = \pi_{t_id^{-1}}(g)$, which is conjugate to $\pi_{s_i}(g)$ if $t_id^{-1} \langle c \rangle = s_i \langle c \rangle$. This does not have to be the case for the set $T_c = \{s_i\}_{i\in I_1}$ computed in Case~2, but we know that for each $i\in I_1 \cup I_2$ there is a unique $k\in I_1 \cup I_2$ such that $t_id^{-1} \langle c \rangle = s_k \langle c \rangle$. Hence, for each $i\in I_1$, it is enough to check for all $k\in I_1$ whether 
$$\pi_{t_i}(f) \text{ and } \pi_{s_k}(g) \text{ are conjugate. }$$ 
If for each $i\in I_1$ there is some $k\in I_1$ for which this is true, then $x$ and $y$ are conjugate. Otherwise, they are not. Note that the above computations amount to solving $L^2$ instances of the Conjugacy Problem in $A$ and so determining whether $x$ and $y$ are conjugate can be done in time 

\begin{equation}\label{eqn: complexity case 3.1}
O\big( L^2 T_{PB}(L) + L^2T_{CA}(L^2) \big).
\end{equation}

\item[2)] \emph{The order of $b$ is infinite.} Let $k$ be a
fixed integer such that $\pi_{t_k}(f) \neq 1$ (such a $k$ must be found
already in the beginning of Case $3$). We proceed to check that
$\pi_{t_k}(f) = \pi_{t_k}^{(d)}(g)$ without finding $d$. Assume that $\pi_{t_k}^{(d)}(g) = 1$ as otherwise, by Theorem~\ref{thm: conjugacy criterion for wreath}, we can conclude that $x$ and $y$ are not conjugate. Since $\pi_{t_k}^{(d)}(g) = \prod_j g(t_kb^jd^{-1}) \neq 1$, there is some integer $l$ for which $g(t_kb^ld^{-1}) \neq 1$. Then $t_kb^ld^{-1} =
\beta_p$ for some $\beta_p \in \supp(g)$ and so $d = \beta_p^{-1}
t_k b^l$. It would suffice to check for all $d$ of the form
$d=\beta_p^{-1} t_k b^l$ such that $db = cd$ whether $\pi_{t_i}(f) =
\pi_{t_i}^{(d)}(g)$.

In order to check the former, we need to check for all $\beta_p \in
\supp(g)$ whether $\beta_p^{-1} t_k b^{l}b=c\beta_p^{-1} t_k b^l$,
i.e., it is enough to check whether $\beta_p^{-1} t_k b =
c\beta_p^{-1} t_k$. These are $m$ instances of the Word Problem in $B$ which do not involve $l$, so they can be decided in time $mT_{CB}(6L)$. Thus checking whether $d$ satisfies $db = cd$ can be done in time $O\big(L T_{CB}(L) \big)$.

It remains to check whether $\pi_{t_i}(f) = \pi_{t_i}^{(d)}(g)$. Notice that
\begin{eqnarray*}
\pi_{t_i}^{(d)}(g) &=& \prod\limits_{j=-\infty}^{\infty}g(t_i
b^jd^{-1})\phantom{t_k^{-1}\beta_p} =
\prod\limits_{j=-\infty}^{\infty} g(t_i b^j b^{-l}t_k^{-1}\beta_p)
\\ &=& \prod\limits_{j=-\infty}^{\infty} g(t_i b^{j-l}
t_k^{-1}\beta_p) \phantom{_k} = \prod\limits_{j=-\infty}^{\infty}
g(t_i b^j t_k^{-1}\beta_p) \; = \; \pi_{t_i}^{(\beta_p^{-1}t_k)}(g).
\end{eqnarray*}
So we need to check whether $\pi_{t_i}^{(\beta_p^{-1}t_k)}(g) =
\pi_{t_i}(f)$. Using \ref{subsec: compute pi_i} this can be done in time
\begin{equation}\label{eqn: complexity case 3.2}
O\big( L T_{CB}(L) + T_{CA}(L^2) + LT_{PB}(L)\big).
\end{equation}

\end{description}
The complexity of the conjugacy problem in $A \wr B$ is
$$O\big( L^2T_{CA}(L^2) + L T_{CB}(L) + L^2T_{PB}(L) \big), $$
which is clearly polynomial since $T_{CA}$, $T_{CB}$ and $T_{PB}$
are polynomial.

\end{proof}

\begin{remark} The algorithm described above differs from the algorithm described in \cite{Matthews:1966} in item $2)$ of Case $3$. The original algorithm is not polynomial in this part.
\end{remark}

\smallskip

\section{Complexity of the Conjugacy Search Problem in Wreath Products}
\label{sec:algorithm wreath CSP}

We use the same notation as in the previous section. The following result is a corollary of several propositions in \cite{Matthews:1966}, together with their proofs. 

\begin{lemma}
\label{lemma: constructing z=dh} Let $A$ and $B$ be finitely
generated groups and let $x=bf$, $y=cg$ be conjugate in $A\wr B$.
Then $z=dh\in A\wr B$ conjugates $x$ to $y$ if and only if $z$ satisfies
\begin{enumerate}
\item $db=cd$ in $B$;
\item when the order of $b$ is finite, $h$ satisfies
\begin{equation}
\label{eqn: h when b is finite order}
h(t_ib^k) = \left(\prod\limits_{j=0}^{k}
g(t_ib^jd^{-1})\right)^{-1}\alpha_i\prod\limits_{j=0}^{k}f(t_ib^j),
\end{equation}
where $\alpha_i$ is such that $\pi_{t_i}^{(d)}(g) =
\alpha_i\pi_{t_i}(f)\alpha_i^{-1}$;
\item when the order of $b$ is infinite, $h$ satisfies
\begin{equation}
\label{eqn: h when b is infinite order}
h(t_ib^k) = \left(\prod\limits_{j=0}^{k}
g(t_ib^jd^{-1})\right)^{-1}\prod\limits_{j=0}^{k}f(t_ib^j).
\end{equation}
\end{enumerate}
\end{lemma}

Note that it follows from \cite{Matthews:1966} that the formulas (\ref{eqn: h when b is finite order}) and (\ref{eqn: h when b is infinite order}) define $h(\beta)$ for all $\beta \in B$ and do not depend on the choice of coset representatives. With this, we can now prove the following theorem.

\begin{theorem}\label{thm: CSP in wreath products is poly}
Let $A$ and $B$ be finitely generated groups such that the following
hold:
 \begin{enumerate}
 \item[1)] there are algorithms which solve the Conjugacy Search Problem in $A$ and in $B$ with polynomial time functions,
 $T_{CSA}$, $T_{CSB}$, respectively;
 \item[2)] there is an algorithm with polynomial time function $T_{PB}$ for the Power Problem in $B$.
 \end{enumerate}
 Then the Conjugacy Search Problem in $A\wr B$ is solvable with complexity
$$O(T_{CSB}(L) + T_{CSA}(L)+ L^2T_{PB}(L)), $$
where $L=|x| +|y|$ is the length of the input pair $x,y \in A\wr B$.
\end{theorem}

\begin{proof}

Let $x=bf$, $y=cg$ be conjugate in $A\wr B$ (this can be checked in polynomial time using Theorem~\ref{thm: CP in wreath products is poly}). Using the algorithm to solve the Conjugacy Search Problem in $B$, one can find $d\in B$ such that $db=cd$ in time $T_{CSB}(L)$. It remains to show that the function $h$ as described in Lemma~\ref{lemma:constructing z=dh} can be described by a finite set of pairs $\{(b_i, h(b_i))\}$. 

First, assume that the order of $b$ in $B$ is infinite. Let
$$M = \max\{M_{i} \mid t_ib^{M_{i}} \in \supp(f)\cup\supp(g), \text{ and } i\in I_1 \}.$$
We show that $M$ can be found in polynomial time. For each $b_j\in \supp(f)\cup \supp(g)$ and for each $t_i \in T_b$, compute $M_{ij}$ such that $t_ib^{M_{ij}} = b_j$. This can be done in time $O(L^2T_{PB}(L))$. Let $M = \max\{M_i \mid b_j \in \supp(f)\cup \supp(g)\}$. Then $M = \max\{ M_i \mid i\in I_1\}$ can be computed in $O(L^2T_{PB}(L))$ steps. Consider the following cases.

\begin{enumerate}
\item $k \geq M$. Then $h(t_ib^k) = \big(\pi_{t_i}^{(d)}(g)\big)^{-1} \pi_{t_i}(f)
=1$, by Theorem~\ref{thm: conjugacy criterion for wreath}. Hence $h(t_ib^k) = 1$.
\item $k<M$.
\begin{enumerate}
\item If $t_i \notin \supp(f)\cup \supp(g)$ and $t_id^{-1} \notin \supp(f)\cup \supp(g)$, then $f(t_ib^j) =1$ and $g(t_id^{-1}b^j) = 1$ for all $j$ and hence $h(\tilde{t_i}b^k)=1$.
\item If $t_i\in \supp(f)\cup \supp(g)$, but $t_id^{-1}\notin \supp(f)\cup \supp(g)$, then
$$h(\tilde{t_i}b^k) = \left( \prod\limits_{j\leq k} g(\tilde{t_i}d^{-1} c^j) \right)^{-1}
\prod\limits_{j\leq k} f(t_ib^j) = \prod\limits_{j\leq k} f(t_ib^j),$$ 
which can be computed in time $O\big(MLT_{PB}(L))\big)$.
\item If $t_i \notin \supp(f)\cup \supp(g)$, but $t_id^{-1}\in \supp(f)\cup \supp(g)$, $h(\tilde{t_i}b^k) = \prod\limits_{j\leq k} g(\tilde{t_i}b^jd^{-1})$ which can be similarly computed in time $O\big(MLT_{PB}(L))\big)$.
\item If $t_i, t_id^{-1}\in \supp(f)\cup \supp(g)$, then 
$$h(\tilde{t_i}b^k) = \left(\prod\limits_{j \leq k} g(t_ib^jd^{-1})\right)^{-1} \prod\limits_{j\leq k} f(t_ib^j)$$
can be computed in time $O\big(MLT_{PB}(L)\big)$.
\end{enumerate}
Thus, if $k< M$, $h(t_ib^k)$ can be computed in time $O\big(MLT_{PB}(L)\big)$. It is clear from the definition of $M$ that $M < L$, so one can compute $h(t_ib^k)$ in time $O(L^2T_{PB}(L)). $
\end{enumerate}

Assume that the order of $b$ is finite, say $N$. Using the algorithm to solve the Conjugacy Search Problem in $A$, one can find in time $T_{CSA}(L^2)$, for each $i\in I_1$, an $\alpha_i \in A$ such that $\pi_{t_i}^{(d)}(g) = \alpha_i \pi_{t_i}(f)\alpha_i^{-1}$. Then $h(t_ib^k) = \left(\prod\limits_{j=0}^{k} g(t_ib^jd^{-1})\right)^{-1} \alpha_i \prod\limits_{j=0}^{k}f(t_ib^j)$ can be found in time $O(T_{CSA}(L) + L^2T_{PB}(L))$ by arguing as in the infinite order case (here instead of $M$, we use the order $N$ of $b$). 

Thus the conjugacy search problem in $A\wr B$ is solvable in time $$O\big( T_{CSB}(L) + T_{CSA}(L) + L^2T_{PB}(L) \big).$$

\end{proof}

\smallskip

\section{Complexity of the Conjugacy and Conjugacy Search Problems in Free Solvable Groups}
\label{sec:algorithm free solvable}

By Corollary~\ref{cor: reduction to CP in wreath is poly} the Conjugacy Problem in free solvable groups can be reduced in polynomial time to the Conjugacy Problem in a wreath product. Then the result from Section~\ref{sec:algorithm wreath} can be applied to deduce that the Conjugacy Problem in free solvable groups is solvable in polynomial time. Though the bound for the complexity will be polynomial, the degree of the polynomial will depend on the degree of solvability (this is because of the factor of $L$ in front of $T_{CB}(L)$ in (\ref{eqn: complexity wreath prod})). However, by making a modification to the algorithm, the complexity of the Conjugacy Problem in free solvable groups is shown to be a polynomial of degree eight.

\begin{theorem}\label{thm: CP in wreath products is poly - modified}
The Conjugacy Problem in a wreath product $A\wr B$, in which $A$ is
abelian is in
$$O\big( T_{CA}(L^2) + T_{CB}(L) + L^2T_{PB}(L) \big),$$
where $L$ is the length of the input pair $(x,y)$.
\end{theorem}

\begin{proof}
The algorithm is similar to the one in Theorem~\ref{thm: CP in wreath products is poly}. The only alteration to be made is in Case~$3$, where the order of $b$ is infinite. Let  $\{s_i\}_{i\in I_1}$ be the set of coset representatives computed in Case~$2$. Then $\pi_{s_i}(g)$ is conjugate to $\pi_{t_i}^{(d)}(g)$. Since $A$ is now abelian, $\pi_{s_i}(g) = \pi_{t_i}^{(d)}(g)$. Thus $\pi_{t_i}(f) = \pi_{t_i}^{(d)}(g)$ if and only if $\pi_{t_i}(f) = \pi_{s_i}(g)$. Checking this requires

\begin{equation}\label{eqn: complexity case 3.2 modified}
O \big( |x|T_{PB}(|x|) + |y|T_{PB}(|y|) + T_{CA}(|x|^2 + |y|^2)
\big).
\end{equation}
As a result the overall complexity of the modified algorithm is
\begin{equation*}
O\big( T_{CA}(L^2) + T_{CB}(L) + L^2T_{PB}(L) \big).
\end{equation*}

\end{proof}

\begin{theorem}
The Conjugacy Problem in $S_{d,r}$ is in $O \big( rdL^8 \big)$,
where $L = |x|+|y|$ is the input length.
\end{theorem}

\begin{proof}
We proceed by induction on the degree of solvability, $d$. The base
case is the abelian group $F/F'$, where the Conjugacy Problem is in $O(rL)$. Now suppose there is an algorithm, which solves the Conjugacy Problem in $F/F^{(d)}$ in $O\big(rd L^8 \big)$. By Corollary~\ref{cor: reduction to CP in wreath is poly}, one can reduce the Conjugacy Problem in $F/F^{(d+1)}$ to the Conjugacy Problem in $F/F' \wr F/F^{(d)}$ in time $O(rd L^3)$. Since $F/F'$ is abelian, we apply Theorem~\ref{thm: CP in wreath products is poly - modified}. In order to do this we need polynomial bounds for the Conjugacy Problems of $F/F'$, $F/F^{(d)}$
and the Power Problem in $F/F^{(d)}$.

The Conjugacy Problem in $F/F'$ is in $O(rL)$. By the induction hypothesis, there is an algorithm which solves the Conjugacy Problem in $F/F^{(d)}$ in $O\big(rdL^8 \big)$. By
Theorem~\ref{prop:power problem is O(dr w^3)} there is an algorithm which solves the Power Problem in $F/F^{(d)}$ in $O(rdL^6)$. Then from Theorem~\ref{thm: CP in wreath products is poly - modified}, the complexity of the Conjugacy Problem in $F/F^{(d+1)}$ is

\begin{equation*}
O \big( rL^2 + rdL^8 + L^2rdL^6 \big),
\end{equation*}

It is easily seen now that the complexity of the Conjugacy Problem in free solvable groups is

\begin{equation*}
O \big( rdL^8 \big).
\end{equation*}

\end{proof}

Since all the proofs of the decidability results are constructive
one can also deduce the following theorem.

\begin{theorem}
The Conjugacy Search Problem in $S_{d,r}$ is solvable in time $O\big( rdL^8 \big)$, where $L = |x|+|y|$ is the input length.
\end{theorem}

\begin{proof}
Again we proceed by induction on the degree of solvability $d$, this
time making sure that at each step we are effectively finding the
required object. When $d=1$ the group is abelian and so deciding the
Conjugacy Search Problem there is trivial -- two words are conjugate
if and only if the identity is a conjugator. Now suppose that there
is an algorithm running in time $O(rd L^8)$, which, if two words $\bar{x}, \bar{y} \in F/F^{(d)}$ are conjugate, exhibits a conjugator. We proceed to describe an algorithm which
does the same for two conjugate elements $x,y \in F/F^{(d+1)}$ given as products of generators of $F$. As before, by Corollary~\ref{cor: reduction to CP in wreath is poly}, we reduce the Conjugacy Problem in $F/F^{(d+1)}$ to the Conjugacy Problem in $F/F^{\prime} \wr F/F^{(d)}$. Hence by Theorem~\ref{thm: CSP in wreath products is poly} there is an algorithm running in time $O(rdL^8)$, which finds a conjugator for $\varphi(x)$ and $\varphi(y)$. The proof of Theorem~2 in \cite{Remeslennikov-Sokolov:1970} gives a pre-image $s \in F/F^{(d+1)}$ for this conjugator. One can see easily that computing $s$ can be done in time $O(r(d+1)L^3)$. Thus, the overall complexity of this algorithm is

$$O\big( r(d+1)L^8 \big). $$
\end{proof}

\smallskip
\smallskip
\smallskip

\bibliography{arxivCP_free_metabelian}

\end{document}